\documentclass[12pt]{amsart}
\usepackage{amssymb}
\begin{document}
\theoremstyle{plain}
\newtheorem{thm}{Theorem}[section]
\newtheorem{theorem}[thm]{Theorem}
\newtheorem{lemma}[thm]{Lemma}
\newtheorem{corollary}[thm]{Corollary}
\newtheorem{proposition}[thm]{Proposition}
\theoremstyle{definition}
\newtheorem{notations}[thm]{Notations}
\newtheorem{remark}[thm]{Remark}
\newtheorem{remarks}[thm]{Remarks}
\newtheorem{definition}[thm]{Definition}
\newtheorem{claim}[thm]{Claim}
\newtheorem{assumption}[thm]{Assumption}
\numberwithin{equation}{thm}
\newcommand{\zar}{{\rm zar}}
\newcommand{\an}{{\rm an}}
\newcommand{\red}{{\rm red}}
\newcommand{\codim}{{\rm codim}}
\newcommand{\rank}{{\rm rank}}
\newcommand{\Pic}{{\rm Pic}}
\newcommand{\Div}{{\rm Div}}
\newcommand{\Hom}{{\rm Hom}}
\newcommand{\im}{{\rm Im}}
\newcommand{\Spec}{{\rm Spec}}
\newcommand{\sing}{{\rm sing}}
\newcommand{\reg}{{\rm reg}}
\newcommand{\Char}{{\rm char}}
\newcommand{\Tr}{{\rm Tr}}
\newcommand{\Gal}{{\rm Gal}}
\newcommand{\Min}{{\rm Min \ }}
\newcommand{\Max}{{\rm Max \ }}
\newcommand{\soplus}[1]{\stackrel{#1}{\oplus}}
\newcommand{\dlog}{{\rm dlog}\,}    
\newcommand{\limdir}[1]{{\displaystyle{\mathop{\rm
lim}_{\buildrel\longrightarrow\over{#1}}}}\,}
\newcommand{\liminv}[1]{{\displaystyle{\mathop{\rm
lim}_{\buildrel\longleftarrow\over{#1}}}}\,}
\newcommand{\boxtensor}{{\Box\kern-9.03pt\raise1.42pt\hbox{$\times$}}}
\newcommand{\sext}{\mbox{${\mathcal E}xt\,$}}
\newcommand{\shom}{\mbox{${\mathcal H}om\,$}}
\newcommand{\coker}{{\rm coker}\,}
\renewcommand{\iff}{\mbox{ $\Longleftrightarrow$ }}
\newcommand{\onto}[1]{\mbox{$\,\>#1 >>\hspace{-.5cm}\to\hspace{.15cm}$}}
\catcode`\@=11
\def\opn#1#2{\def#1{\mathop{\kern0pt\fam0#2}\nolimits}}
\def\bold#1{{\bf #1}}%
\def\underrightarrow{\mathpalette\underrightarrow@}
\def\underrightarrow@#1#2{\vtop{\ialign{$##$\cr
 \hfil#1#2\hfil\cr\noalign{\nointerlineskip}%
 #1{-}\mkern-6mu\cleaders\hbox{$#1\mkern-2mu{-}\mkern-2mu$}\hfill
 \mkern-6mu{\to}\cr}}}
\let\underarrow\underrightarrow
\def\underleftarrow{\mathpalette\underleftarrow@}
\def\underleftarrow@#1#2{\vtop{\ialign{$##$\cr
 \hfil#1#2\hfil\cr\noalign{\nointerlineskip}#1{\leftarrow}\mkern-6mu
 \cleaders\hbox{$#1\mkern-2mu{-}\mkern-2mu$}\hfill
 \mkern-6mu{-}\cr}}}
\let\amp@rs@nd@\relax
\newdimen\ex@
\ex@.2326ex
\newdimen\bigaw@
\newdimen\minaw@
\minaw@16.08739\ex@
\newdimen\minCDaw@
\minCDaw@2.5pc
\newif\ifCD@
\def\minCDarrowwidth#1{\minCDaw@#1}
\newenvironment{CD}{\@CD}{\@endCD}
\def\@CD{\def\A##1A##2A{\llap{$\vcenter{\hbox
 {$\scriptstyle##1$}}$}\Big\uparrow\rlap{$\vcenter{\hbox{%
$\scriptstyle##2$}}$}&&}%
\def\V##1V##2V{\llap{$\vcenter{\hbox
 {$\scriptstyle##1$}}$}\Big\downarrow\rlap{$\vcenter{\hbox{%
$\scriptstyle##2$}}$}&&}%
\def\={&\hskip.5em\mathrel
 {\vbox{\hrule width\minCDaw@\vskip3\ex@\hrule width
 \minCDaw@}}\hskip.5em&}%
\def\verteq{\Big\Vert&&}%
\def\noarr{&&}%
\def\vspace##1{\noalign{\vskip##1\relax}}\relax\let\amp@rs@nd@&\iffalse}\fi
 \CD@true\vcenter\bgroup\relax\let\\=\cr\iffalse}\fi\tabskip\z@skip\baselineskip20\ex@
 \lineskip3\ex@\lineskiplimit3\ex@\halign\bgroup
 &\hfill$\m@th##$\hfill\cr}
\def\@endCD{\cr\egroup\egroup}
\def\>#1>#2>{\amp@rs@nd@\setbox\z@\hbox{$\scriptstyle
 \;{#1}\;\;$}\setbox\@ne\hbox{$\scriptstyle\;{#2}\;\;$}\setbox\tw@
 \hbox{$#2$}\ifCD@
 \global\bigaw@\minCDaw@\else\global\bigaw@\minaw@\fi
 \ifdim\wd\z@>\bigaw@\global\bigaw@\wd\z@\fi
 \ifdim\wd\@ne>\bigaw@\global\bigaw@\wd\@ne\fi
 \ifCD@\hskip.5em\fi
 \ifdim\wd\tw@>\z@
 \mathrel{\mathop{\hbox to\bigaw@{\rightarrowfill}}\limits^{#1}_{#2}}\else
 \mathrel{\mathop{\hbox to\bigaw@{\rightarrowfill}}\limits^{#1}}\fi
 \ifCD@\hskip.5em\fi\amp@rs@nd@}
\def\<#1<#2<{\amp@rs@nd@\setbox\z@\hbox{$\scriptstyle
 \;\;{#1}\;$}\setbox\@ne\hbox{$\scriptstyle\;\;{#2}\;$}\setbox\tw@
 \hbox{$#2$}\ifCD@
 \global\bigaw@\minCDaw@\else\global\bigaw@\minaw@\fi
 \ifdim\wd\z@>\bigaw@\global\bigaw@\wd\z@\fi
 \ifdim\wd\@ne>\bigaw@\global\bigaw@\wd\@ne\fi
 \ifCD@\hskip.5em\fi
 \ifdim\wd\tw@>\z@
 \mathrel{\mathop{\hbox to\bigaw@{\leftarrowfill}}\limits^{#1}_{#2}}\else
 \mathrel{\mathop{\hbox to\bigaw@{\leftarrowfill}}\limits^{#1}}\fi
 \ifCD@\hskip.5em\fi\amp@rs@nd@}
\newenvironment{CDS}{\@CDS}{\@endCDS}
\def\@CDS{\def\A##1A##2A{\llap{$\vcenter{\hbox
 {$\scriptstyle##1$}}$}\Big\uparrow\rlap{$\vcenter{\hbox{%
$\scriptstyle##2$}}$}&}%
\def\V##1V##2V{\llap{$\vcenter{\hbox
 {$\scriptstyle##1$}}$}\Big\downarrow\rlap{$\vcenter{\hbox{%
$\scriptstyle##2$}}$}&}%
\def\={&\hskip.5em\mathrel
 {\vbox{\hrule width\minCDaw@\vskip3\ex@\hrule width
 \minCDaw@}}\hskip.5em&}
\def\verteq{\Big\Vert&}
\def\novarr{&}
\def\noharr{&&}
\def\SE##1E##2E{\slantedarrow(0,18)(4,-3){##1}{##2}&}
\def\SW##1W##2W{\slantedarrow(24,18)(-4,-3){##1}{##2}&}
\def\NE##1E##2E{\slantedarrow(0,0)(4,3){##1}{##2}&}
\def\NW##1W##2W{\slantedarrow(24,0)(-4,3){##1}{##2}&}
\def\slantedarrow(##1)(##2)##3##4{%
\thinlines\unitlength1pt\lower 6.5pt\hbox{\begin{picture}(24,18)%
\put(##1){\vector(##2){24}}%
\put(0,8){$\scriptstyle##3$}%
\put(20,8){$\scriptstyle##4$}%
\end{picture}}}
\def\vspace##1{\noalign{\vskip##1\relax}}\relax\let\amp@rs@nd@&\iffalse}\fi
 \CD@true\vcenter\bgroup\relax\let\\=\cr\iffalse}\fi\tabskip\z@skip\baselineskip20\ex@
 \lineskip3\ex@\lineskiplimit3\ex@\halign\bgroup
 &\hfill$\m@th##$\hfill\cr}
\def\@endCDS{\cr\egroup\egroup}
\newdimen\TriCDarrw@
\newif\ifTriV@
\newenvironment{TriCDV}{\@TriCDV}{\@endTriCD}
\newenvironment{TriCDA}{\@TriCDA}{\@endTriCD}
\def\@TriCDV{\TriV@true\def\TriCDpos@{6}\@TriCD}
\def\@TriCDA{\TriV@false\def\TriCDpos@{10}\@TriCD}
\def\@TriCD#1#2#3#4#5#6{%
\setbox0\hbox{$\ifTriV@#6\else#1\fi$}
\TriCDarrw@=\wd0 \advance\TriCDarrw@ 24pt
\advance\TriCDarrw@ -1em
\def\SE##1E##2E{\slantedarrow(0,18)(2,-3){##1}{##2}&}
\def\SW##1W##2W{\slantedarrow(12,18)(-2,-3){##1}{##2}&}
\def\NE##1E##2E{\slantedarrow(0,0)(2,3){##1}{##2}&}
\def\NW##1W##2W{\slantedarrow(12,0)(-2,3){##1}{##2}&}
\def\slantedarrow(##1)(##2)##3##4{\thinlines\unitlength1pt
\lower 6.5pt\hbox{\begin{picture}(12,18)%
\put(##1){\vector(##2){12}}%
\put(-4,\TriCDpos@){$\scriptstyle##3$}%
\put(12,\TriCDpos@){$\scriptstyle##4$}%
\end{picture}}}
\def\={\mathrel {\vbox{\hrule
   width\TriCDarrw@\vskip3\ex@\hrule width
   \TriCDarrw@}}}
\def\>##1>>{\setbox\z@\hbox{$\scriptstyle
 \;{##1}\;\;$}\global\bigaw@\TriCDarrw@
 \ifdim\wd\z@>\bigaw@\global\bigaw@\wd\z@\fi
 \hskip.5em
 \mathrel{\mathop{\hbox to \TriCDarrw@
{\rightarrowfill}}\limits^{##1}}
 \hskip.5em}
\def\<##1<<{\setbox\z@\hbox{$\scriptstyle
 \;{##1}\;\;$}\global\bigaw@\TriCDarrw@
 \ifdim\wd\z@>\bigaw@\global\bigaw@\wd\z@\fi
 \mathrel{\mathop{\hbox to\bigaw@{\leftarrowfill}}\limits^{##1}}
 }
 \CD@true\vcenter\bgroup\relax\let\\=\cr\iffalse}\fi
 \tabskip\z@skip\baselineskip20\ex@
 \lineskip3\ex@\lineskiplimit3\ex@
 \ifTriV@
 \halign\bgroup
 &\hfill$\m@th##$\hfill\cr
#1&\multispan3\hfill$#2$\hfill&#3\\
&#4&#5\\
&&#6\cr\egroup%
\else
 \halign\bgroup
 &\hfill$\m@th##$\hfill\cr
&&#1\\%
&#2&#3\\
#4&\multispan3\hfill$#5$\hfill&#6\cr\egroup
\fi}
\def\@endTriCD{\egroup}
\newcommand{\sA}{{\mathcal A}}
\newcommand{\sB}{{\mathcal B}}
\newcommand{\sC}{{\mathcal C}}
\newcommand{\sD}{{\mathcal D}}
\newcommand{\sE}{{\mathcal E}}
\newcommand{\sF}{{\mathcal F}}
\newcommand{\sG}{{\mathcal G}}
\newcommand{\sH}{{\mathcal H}}
\newcommand{\sI}{{\mathcal I}}
\newcommand{\sJ}{{\mathcal J}}
\newcommand{\sK}{{\mathcal K}}
\newcommand{\sL}{{\mathcal L}}
\newcommand{\sM}{{\mathcal M}}
\newcommand{\sN}{{\mathcal N}}
\newcommand{\sO}{{\mathcal O}}
\newcommand{\sP}{{\mathcal P}}
\newcommand{\sQ}{{\mathcal Q}}
\newcommand{\sR}{{\mathcal R}}
\newcommand{\sS}{{\mathcal S}}
\newcommand{\sT}{{\mathcal T}}
\newcommand{\sU}{{\mathcal U}}
\newcommand{\sV}{{\mathcal V}}
\newcommand{\sW}{{\mathcal W}}
\newcommand{\sX}{{\mathcal X}}
\newcommand{\sY}{{\mathcal Y}}
\newcommand{\sZ}{{\mathcal Z}}
\newcommand{\A}{{\mathbb A}}
\newcommand{\B}{{\mathbb B}}
\newcommand{\C}{{\mathbb C}}
\newcommand{\D}{{\mathbb D}}
\newcommand{\E}{{\mathbb E}}
\newcommand{\F}{{\mathbb F}}
\newcommand{\G}{{\mathbb G}}
\newcommand{\HH}{{\mathbb H}}
\newcommand{\I}{{\mathbb I}}
\newcommand{\J}{{\mathbb J}}
\newcommand{\M}{{\mathbb M}}
\newcommand{\N}{{\mathbb N}}
\renewcommand{\P}{{\mathbb P}}
\newcommand{\Q}{{\mathbb Q}}
\newcommand{\R}{{\mathbb R}}
\newcommand{\T}{{\mathbb T}}
\newcommand{\U}{{\mathbb U}}
\newcommand{\V}{{\mathbb V}}
\newcommand{\W}{{\mathbb W}}
\newcommand{\X}{{\mathbb X}}
\newcommand{\Y}{{\mathbb Y}}
\newcommand{\Z}{{\mathbb Z}}
\title[Families of projective
manifolds over curves]{On the isotriviality of families
of projective manifolds over curves}
\author{Eckart Viehweg}
\address{Universit\"at GH Essen, FB6 Mathematik, 45117 Essen, Germany}
\email{ viehweg@uni-essen.de}
\thanks{This work has been partly supported by the DFG
Forschergruppe ``Arithmetik und Geometrie''}
\author[Kang Zuo]{Kang Zuo${}^*$}
\address{Universit\"at Kaiserslautern, FB Mathematik, 67653 Kaiserslautern, Germany}
\email{zuo@mathematik.uni-kl.de}
\thanks{${}^*$ Supported by a Heisenberg-fellowship, DFG}
\maketitle
Let $Y$ be a non-singular curve, $X$ a manifold, both projective
and defined over $\C$, and let $f: X \to Y$ be a surjective
morphism with connected general fibre $F$. We fix a reduced
divisor $S$ on $Y$ which contains the discriminant divisor of $f$, i.e.
a reduced divisor with
$$
f_0 = f|_{X_0} : X_0 \>>> Y_0
$$
smooth, for $Y_0 = Y - S$ and $X_0 = f^{-1} (Y_0)$.
Recall that $f$ is birationally isotrivial, if $X \times_Y {\rm Spec} \
\overline{\C(Y)}$ is birational to $F \times {\rm Spec}
\overline{\C(Y)}$.

\begin{theorem} \label{main1}
Assume that $f$ is not birationally isotrivial, and that
one of the following conditions holds true:
\begin{enumerate}
\item[a)] $\kappa (F) = \dim (F)$.
\item[b)] $F$ has a minimal model $F'$ with $K_{F'}$ semi-ample.
\end{enumerate}
Then $f$ has at least
\begin{enumerate}
\item[i)] three singular fibres if $Y=\P^1$.
\item[ii)] one singular fibre if $Y$ is an elliptic curve.
\end{enumerate}
\end{theorem}

If $\kappa (F) = \dim (F) = 2$ or if $\omega_F$ is ample, part
ii) of \ref{main1} has been shown by L. Migliorini \cite{Mig} and part
i) is due to S. Kov\'acs \cite{Kov1}, \cite{Kov2}. If $\dim (F) = 2$ and
$\kappa (F) \leq 1$, both i) and ii) have recently been shown
by K. Oguiso and the first named author \cite{O-V}. The proof is based
on the observation that the additional singular fibres,
showing up in certain cyclic coverings of $X$, can be neglected
in vanishing theorems.

The same principle is exploited in the proof of \ref{main1},
replacing global vanishing theorems by the negativity of the
kernel of the Kodaira-Spencer map,
as used by J. Jost and the second named author in
\cite{J-Z} (see also \cite{Pet}).

Since the total space $X$ of a smooth isotrivial morphism to an
elliptic curve has at most Kodaira dimension $\dim(X)-1$, part
ii) in \ref{main1} implies that there are no smooth morphisms from a
manifold of general type to an elliptic curve. A similar
argument shows that a surjective morphism from a manifold of
general type to $\P^1$ must have at least three singular fibres,
an affirmative answer to a question posed by F. Catanese and M.
Schneider. In fact, one has a stronger result.
The assumptions a) or b) in \ref{main1} are just needed to guarantee
that for some $\nu\gg 0$ the determinant of $f_*\omega_{X/Y}^\nu$
is ample. For $Y=\P^1$ this necessarily holds true, if $X$ has a
non-negative Kodaira dimension.

\begin{theorem}\label{main3}
Let $X$ be a complex projective manifold of non-negative Kodaira
dimension. Then a surjective morphism $f:X \to \P^1$
has at least 3 singular fibres.
\end{theorem}

If $X$ is a curve, this is just the wellknown fact
that a morphism to $\P^1$ from a curve of genus larger than or equal to
one has to ramify over at least three points.
For surfaces \ref{main3} follows from the bounds of A.
Parshin \cite{Par} and A. Arakelov \cite{Ara}. For
threefolds of general type, or if $\omega_X$ is ample,
\ref{main3} is an easy corollary of the results in \cite{Kov2} (see
also \cite{B-V}). For threefolds of lower Kodaira dimension one
can use \cite{O-V} instead.

As a byproduct, the proof of \ref{main1} gives some explicit
bounds, generalizing the ones obtained by A. Parshin \cite{Par}
and A. Arakelov \cite{Ara} for families of curves, by E. Bedulev \cite{Bed},
for $F$ an elliptic surface, and by E. Bedulev and the first author
\cite{B-V} under the assumption that $F$ is
canonically polarized, or a surface of general type.
We write $s = \deg(S)$, and $n = \dim( F)$. The genus of $Y$
is denoted by $g$, and $\delta$ is the number of singular
fibres of $f$, which are not semistable, i.e. not reduced
normal crossing divisors.

\begin{theorem} \label{main2}
Assume that $f$ is not birationally isotrivial, that $\omega_F$ is
semi-ample, and that $h(t)$ is the Hilbert polynomial for a
polarization of $F$. Then there exist constants $\nu$ and $e$,
depending only on $h$, with
$$
\frac{\deg (f_* \omega^{\nu}_{X/Y})}{\rank (f_*
\omega^{\nu}_{X/Y})} \leq ( n \cdot (2 g - 2 + s) + \delta ) \cdot \nu
\cdot e.
$$
In particular, if $f$ is semistable,
$$
\frac{\deg (f_* \omega^{\nu}_{X/Y})}{\rank (f_*
\omega^{\nu}_{X/Y})} \leq n \cdot (2 g - 2 + s) \cdot \nu \cdot e.
$$
\end{theorem}
As explained in \cite{B-V}, section 4, such bounds imply a
generalization of the Shafarevich conjecture, saying that
for given $Y$ and $S$ there are finitely many deformation types,
if there exists a nice compactification of the
corresponding moduli scheme. Unfortunately such a compactification
has only been constructed for surfaces of general type.

We thank Keiji Oguiso for helpful remarks and comments, and
Osamu Fujino for pointing out some ambiguities in the first
version of section three.
Sandor Kov\'acs informed us, immediately after he got a
copy of this article, that by different methods he obtained
results, partly overlapping with the ones presented here
(see \cite{Kov3}).

\section{Hodge bundles and the Kodaira Spencer map}

Let $Y$ be a complex manifold and let $S$ be a normal crossing
divisor. A variation $\V_0$ of polarized Hodge structures of weight $k$
on $Y_0=Y-S$ gives rise to
$$
E_0 = {\rm gr}_F(\V\otimes \sO_{Y_0}) = \bigoplus_{p+q=k} E^{p,q}_0,
$$
together with a Higgs structure
$\theta_0 = \oplus \theta_{p,q} : E_0 \to E_0\otimes \Omega^1_{Y_0}$.
\begin{lemma}\label{hodge-bundles}
If  $ \sN \subset E_0^{p,q}$ is  a sub-bundle with $\theta_{p,q}(\sN)=0,$
then the curvature of the restricted Hodge metric on $\sN$ is negative
semidefinite.
\end{lemma}
In fact, the negativity of the curvature of the restricted Hodge
metric on $\det(\sN)$, which will be the only case used in this note,
as well as the next lemma \ref{hodge-bundles2}, follow from
\cite{Sim} and can also be found in \cite{J-Z}, Lemma 1. For the
convenience of the reader we sketch the proof.
\begin{proof}
Let $ \Theta(E_0,h)$ denote the curvature form of the Hodge metric $ h$
on $ E_0.$ Then by \cite{Gri}, chapter II, we have
$$
\Theta(E_0,h)+\theta\wedge
\bar\theta_h+\bar\theta_h\wedge\theta=0,
$$
where $\bar\theta_h$ is the complex conjugation of $ \theta$
with respect to $ h.$
$h$ restricts to a metric $ h|_\sN,$ and induces a
$\sC^{\infty}-$decomposition $ E_0=\sN\oplus \sN^{\perp}.$
One obtains
$$ \Theta(\sN,h)=\Theta(E_0,h)|_{\sN}+\bar A_h\wedge A=
-\theta\wedge\bar\theta_h|_{\sN}- \bar
\theta_h\wedge\theta|_{\sN}+\bar A_h\wedge A,
$$
where $ A\in A^{1,0}(\makebox{Hom}(\sN, {\sN}^{\perp}))$ is the
second fundamental form of the sub-bundle $\sN\subset E_0,$ and $
\bar A_h$ is its complex conjugate with respect to $ h.$

Since  $ \theta(\sN)=0,$ we have $ \bar
\theta_h\wedge\theta|_{\sN}=0,$ hence
$$
\Theta(\sN,h)=-\theta\wedge\bar\theta_h|_{\sN}+\bar A_h\wedge A.
$$
$\theta\wedge\bar\theta_h$ is positive semidefinite and $ \bar
A_h \wedge A$ is negative semidefinite, so $ \Theta(\sN,h)$ is
negative semidefinite.
\end{proof}
Suppose the local monodromy of $\V_0$ around the components of
$ S$ are unipotent and let $\sV$ be the Deligne extension of
$\V_0\otimes\sO_{Y_0}$. By \cite{Sch} the F-filtration extends
to a filtration of $\sV$ by subbundles, hence there exists a canonical
extension $ E$ of $ E_0$ to $ Y,$ and $\theta_0$ extends to
$$
\theta=\bigoplus_{p+q=k}\theta_{p,q}: E=\bigoplus_{p+q=k}
E^{p,q} \>>> E\otimes\Omega^1_Y(\log S) =\bigoplus_{p+q=k}
E^{p,q} \otimes \Omega^1_Y(\log S).
$$
\begin{lemma}\label{hodge-bundles2}
Keeping the assumptions made above, suppose $Y$ is a smooth
projective curve. If  $ \sN \subset E^{p,q}$ is a sub-bundle
with $ \theta_{p,q}(\sN)=0,$ then $\deg (\sN)\leq 0.$
\end{lemma}
\begin{proof} Let $\sN^\vee$ be the dual of $\sN$.
We have the projection $ {E^{p,q}}^\vee\to \sN^\vee.$
Note that $ {E^{p,q}}^\vee = {E^{\vee}}^{q,p}$ as a Hodge bundle
of the system of Hodge bundles corresponding to the dual variation of
Hodge structures $\V_0^\vee$. The monodromy of $\V_0^\vee$
around $ S $ is again unipotent. We have the projection
$$ F^{\vee q}\to E^{\vee q,p}\to \sN^\vee, $$
where  $ F^{\vee q}$ is the q-th subbundle in the extended Hodge
filtration of $ \sV^\vee.$

This presentation of $\sN^\vee$ as a quotient of a subbundle of
a variation of Hodge structures allows to apply \cite{Kol}, 5.20.
So the Chern forms of the induced
Hodge metric on $(\sN|_{Y_0})^\vee$ represent the corresponding Chern
classes of $ \sN^\vee.$ From \ref{hodge-bundles} we get in particular
$ \deg ( \sN^\vee)\geq 0$, and hence $ \deg (\sN)\leq 0.$
\end{proof}

Let $g:Z \to Y$ be a surjective morphism between a projective
$n$-dimensional manifold $Z$ and a non-singular curve $Y$,
both defined over the complex numbers. Let $S\subset Y$ be a
divisor such that $g$ is smooth outside of $\Pi=g^{-1}(S)$.
We will assume $\Pi$ to be a normal crossing divisor.
The smooth projective morphism
$$
g_0: Z_0=Z-\Pi \>>> Y-S
$$
obtained by restricting $g$ gives rise to variations of Hodge
structures $\V_0=R^k{g_0}_*\C_{Z_0}$. As explained in
\cite{Zuc}, p. 423, the primitive decomposition of $\V_0$ allows
to define a polarization on $\V_0$. If the fibres of $g$ are
connected and if $g$ is semistable, i.e. if $\Pi$ is reduced,
the local monodromies around points in $S$ are unipotent.

Using the notations introduced above we find
$$
E^{p,q} = R^qg_*\Omega^p_{Z/Y}(\log \Pi).
$$
$\theta_{p,q}:E^{p,q} \to E^{p-1,q+1}\otimes \Omega^1_Y(\log
S)$, which we will call the Kodaira Spencer map, is the
edge-morphism induced by the tautological exact sequence
\begin{equation} \label{a}
0 \to g^*
\Omega^{1}_{Y} (\log S) \otimes \Omega^{p-1}_{Z/Y} (\log \Pi)
\to \Omega^{p}_{Z} (\log \Pi) \to \Omega^{p}_{Z/Y} (\log \Pi)
\to 0.
\end{equation}
It is given by the cup product with the Kodaira Spencer
class, induced by $g$.
\begin{proposition} \label{kernel}
Let $\sN$ be an invertible subsheaf of $E^{p,q}$ with
$\theta_{p,q} (\sN) = 0$. Then
$\deg (\sN) \leq 0.$
\end{proposition}
\begin{proof}
If the fibres of $g$ are connected and if $g$ is semistable,
this is nothing but \ref{hodge-bundles2}.

In general, let $L$ be a finite extension of the
function field $\C(Y)$, containing the Galois hull of the
algebraic closure of $\C(Y)$ in $\C(Z)$, and let
$Y'\to Y$ be the normalization of $Y$ in $L$.
Consider the normalization $\tilde{Z}$ of
$Z\times_YY'$, a desingularization $\varphi':Z' \to \tilde{Z}$
and the induced morphisms
$$
\begin{CDS}
Z' \> \varphi' >> \tilde{Z} \> \tilde{\varphi} >> Z\times_YY'
\> p_1 >> Z \\
\novarr \SE g'\ \ E E \V \tilde{g} V V \SW p_2 W W \novarr \SW W g W \\
\noharr \ Y' \ \> \psi >> Y \ \ \
\end{CDS}
$$
$\varphi=\tilde{\varphi}\circ \varphi'$ and $\psi'= p_1\circ \varphi$.
We will enlarge $S$ such that $Y' \to Y$ is \'etale
over $Y-S$, hence for $S'=\varphi^*S$
$$
\psi^*\Omega^1_Y(\log S) = \Omega^1_{Y'}(\log
S').
$$
If one chooses $L$ large enough, $Z'$ will be the disjoint
union of semistable families over $Y'$, hence $\Pi'={g'}^*S'$ is a reduced
normal crossing divisor, and $\varphi|_{{g'}^{-1}(Y'-S')}$ is
an isomorphism. By the generalized Hurwitz formula \cite{E-V2}, 3.20,
$$
{\psi'}^*\Omega^p_Z(\log \Pi) \subset \Omega^1_{Z'}(\log \Pi'),
$$
and by \cite{E-V1}, Lemme 1.2,
$$
R^q\varphi'_*\Omega^{p}_{Z'}(\log \Pi') =
\left\{ \begin{array}{ll}
\tilde{\varphi}^*p_1^*\Omega^{p}_{Z}(\log \Pi) & \mbox{for \ } q=0 \\
0 & \mbox{for \ } q> 0.
\end{array} \right.
$$
The exact sequence (\ref{a}), for $Z'$ instead of $Z$, and
induction on $p$ allow to show the same for the relative
differential forms, i.e.
$$
R^q\varphi'_*\Omega^{p}_{Z'/Y'}(\log \Pi') =
\left\{ \begin{array}{ll}
\tilde{\varphi}^*p_1^*\Omega^{p}_{Z/Y}(\log \Pi) & \mbox{for \ } q=0 \\
0 & \mbox{for \ } q> 0.
\end{array} \right.
$$
Hence the pullback of the exact sequence (\ref{a}) to
$\tilde{Z}$ is isomorphic to
\begin{multline} \label{a'}
0 \>>> {\varphi'}_* ({g'}^*
\Omega^{1}_{Y'} (\log S') \otimes \Omega^{p-1}_{Z'/Y'} (\log \Pi'))\\
\>>> {\varphi'}_* \Omega^{p}_{Z'} (\log \Pi') \>>>
{\varphi'}_*\Omega^{p}_{Z'/Y'} (\log \Pi')
\>>> 0.
\end{multline}
Writing
$$
{E'}^{p,q} = R^qg'_*\Omega^p_{Z'/Y'}(\log \Pi'),
$$
and $\theta'_{p,q}$ for the edge-morphism, we find
$$
{E'}^{p,q} = R^q{p_2}_*(\tilde{\varphi}_*\sO_{\tilde{Z}}\otimes
p_1^*\Omega^p_{Z/Y}(\log \Pi)).
$$
Moreover, the inclusion $\sO_{Z\times_YY'} \to
\tilde{\varphi}_*\sO_{\tilde{Z}}$ and flat base change give
an inclusion
$$
\psi^*E^{p,q}=R^q{p_2}_*p_1^*\Omega^p_Z(\log \Pi)) \>>> {E'}^{p,q}
$$
and the diagram
$$
\begin{CD}
\psi^*{E}^{p,q} \> \psi^*\theta_{p,q}>> \psi^*{E}^{p-1,q+1}\otimes
\Omega^1_{Y} (\log S)\\
\V\subset VV \V\subset VV\\
{E'}^{p,q} \> \theta'_{p,q}>> {E'}^{p-1,q+1}\otimes
\Omega^1_{Y'} (\log S')
\end{CD}
$$
commutes. In particular, if $\sN$ lies in the kernel of
$\theta_{p,q}$, the sheaf $\psi^*\sN$ lies in the kernel of
$\theta'_{p,q}$. Since we already know \ref{kernel} for
semistable morphisms with connected fibres, we find
$\deg(\sN) \leq 0$.
\end{proof}

\section{Positivity of direct image sheaves}

As in \cite{B-V} or \cite{O-V} a second ingredient in the proof
of \ref{main1} and \ref{main2} will be explicit bounds for the
positivity of certain direct image sheaves.

\begin{definition} \label{nef}
Let $\sE$ be a locally free sheaf and $\sA$ an invertible sheaf
on the curve $Y$.
\begin{enumerate}
\item[a)] $\sE$ is nef if for all finite morphisms $\pi : Z
\to Y$ and all invertible quotients $\sM$ of $\pi^* \sE$ the
degree of $\sM$ is non-negative.
\item[b)] For $\alpha, \beta \in \N - \{ 0 \}$ we write
$$
\sE \succeq \frac{\alpha}{\beta} \sA
$$
if $S^{\beta} (\sE) \otimes \sA^{-\alpha}$ is nef. This is
welldefined, since obviously the latter holds true if and only if
$$
S^{\beta\cdot\mu} (\sE) \otimes \sA^{-\alpha\cdot\mu}
$$
is nef, for some $\mu>0$.
\end{enumerate}
\end{definition}

Nef locally free sheaves on curves,
have already been used in \cite{E-V} to study the height of
points of curves over function fields. In \cite{Vie}, \S 2, and in the
higher dimensional birational classification theory, one
needs positive coherent torsionfree sheaves over higher dimensional
manifolds, and there one often considers weakly positive
sheaves, instead of nef sheaves. In the one-dimensional case,
both notions coincide, and all the properties of weakly
positive sheaves, listed in \cite{Mor} or \cite{Vie} carry over
to nef sheaves on curves. Let us recall one property:

\begin{lemma} \label{twist} Given $d \in \N$, assume that for
all $\mu \in \N$, sufficiently large and divisible,
there exists a covering $\tau : Y'
\to Y$ of degree $\mu$ such that $\tau^* \sE \otimes \sH$ is nef,
for one, hence for all invertible sheaves of degree $d$. Then $\sE$ is nef.
\end{lemma}

\begin{proof}
Let $\pi : Z \to Y$ and $\sM$ be as in \ref{nef}, a), and let
$Z'$ be a component of the normalization of $Z \times_{Y} Y'$.
If
$$
\begin{CD}
Z' \> \tau' >> Z \\
\V \pi' VV \V V \pi V \\
Y' \> \tau >> Y
\end{CD}
$$
are the induced morphisms, then
\begin{multline*}
0 \leq \deg (\tau{'}^{*} \sM \otimes \pi{'}^{*} \sH) = \deg
(\tau') \cdot \deg (\sM) + \deg (\pi')\cdot d\\
\leq \mu \cdot \deg (\sM) + \deg (\pi)\cdot d .
\end{multline*}
This, for all $\mu \in \N - \{ 0 \}$, implies that $\deg (\sM)
\geq 0$.
\end{proof}

Next recall the definition of the (algebraic) multiplier
sheaves. We consider a surjective morphism $f:X\to Y$, with
connected general fibre $F$, where $X$ is an (n+1)-dimensional
complex projective manifold, and $Y$ a non-singular projective
curve. If $\Gamma$ is an effective divisor on $X$,
$$
\omega_{X/Y} \left\{ - \frac{\Gamma}{N} \right\} = \tau_*
\left(\omega_{X'/Y} \left( - \left[ \frac{\Gamma'}{N} \right] \right)\right)
$$
where $\tau: X' \to X$ is any blowing up with $\Gamma' = \tau^*
\Gamma$ a normal crossing divisor (see for example \cite{E-V2}, 7.4,
or \cite{Vie}, section 5.3).

Fujita's positivity theorem (today an easy corollary of
Koll\'ar's vanishing theorem) says that $f_* \omega_{X/Y}$ is
nef. A direct consequence is the following.

\begin{lemma} \label{fujita} Let $\sN$ be an invertible sheaf on
$X$ and $\Gamma$ an effective divisor. Assume that for some $N >
0$ there exists a nef locally free sheaf $\sE$ on $Y$ and a
surjection $f^* \sE \to \sN^N (- \Gamma).$
Then
$$f_* \left(\sN \otimes \omega_{X/Y} \displaystyle \left\{ - \frac{\Gamma}{N}
\right\} \right)$$
is nef.
\end{lemma}

\begin{proof}
Let $p \in Y$ be a point. Then $\sE \otimes \sO_Y (N \cdot p)$
is ample, hence
$$
\sN^N (- \Gamma) \otimes f^* \sO_Y (N \cdot p)
$$
is semi-ample. By \cite{E-V2}, 7.16, the sheaf
$$
f_* \left(\sN \otimes \omega_{X/Y} \left\{ - \frac{\Gamma}{N}
\right\} \right) \otimes \sO_Y (p)
$$
is nef. Since the same holds true over all $Y'$, finite over $Y$
and unramified in $S$, one obtains \ref{fujita} from \ref{twist}
\end{proof}

As an application of \ref{fujita} one obtains, as explained in
\cite{Mor},
\begin{lemma}\label{weak}\ \\[-.4cm]
\begin{enumerate}
\item[i)] $f_* \omega^{\nu}_{X/Y}$ is nef, for all $\nu\geq 0$.
\item[ii)] If $\lambda_\nu= \det(f_* \omega^{\nu}_{X/Y})$
is ample, for some $\nu >1$, then there exists a positive
rational number $\eta$ with $f_* \omega^{\nu}_{X/Y} \succeq \eta \cdot
\lambda_{\nu}$.
\end{enumerate}
\end{lemma}
As in \cite{B-V}, we will need an explicit bound for the
rational number $\eta$ in \ref{weak}, ii). To this aim recall
the following definition, used in \cite{E-V}, \cite{E-V2}, \S \ 7 and
\cite{Vie}, section 5.3.

\begin{definition} \label{e} Let $\sL$ be an invertible sheaf on
$F$ with $H^0 (F, \sL) \neq 0$, and let $\Gamma$ be an effective
divisor. Then
$$
e(\Gamma) = \Min \left\{ N \in \N - \{ 0 \} ; \ \omega_F \left\{-
\frac{\Gamma}{N} \right\} = \omega_F \right\} \ \ \mbox{ \ \ \
\ \ and}
$$
$$
e (\sL) = \Max \left\{ e (\Gamma); \ \Gamma \ \mbox{the zero set of} \
\sigma \in H^0 (F, \sL) - \{ 0 \} \right\}.
$$
\end{definition}

\begin{notations} \label{constants}
For $f: X \to Y$, we choose $\nu > 1$ with $f_*
\omega^{\nu}_{X/Y} = \sE \neq 0$, and a blowing up $\tau: X' \to
X$ such that the fibres of $f' = f \circ \tau$ are normal
crossing divisors, such that
$$
\sL = {\rm Im} (f{'}^* f'_* \omega^{\nu}_{X'/Y} = f{'}^* \sE \>>>
\omega^{\nu}_{X'/Y})
$$
is invertible and $\omega^{\nu}_{X'/Y} = \sL (B)$, for a normal
crossing divisor $B$. Let $F'$ be a general fibre of $f'$. We
define: \\[.2cm]
\indent $e = e (\sL |_{F'})$ \\[.2cm]
\indent $r = \rank ( f'_* \omega^{\nu}_{X'/Y}) = \rank (\sE)$\\[.2cm]
\indent $\lambda = \det (f'_* \omega^{\nu}_{X'/Y}) = \det(\sE)$.
\end{notations}

\begin{proposition} \label{bounds}
If $\displaystyle \lambda$ is ample, $f_* \omega^{\nu}_{X/Y}
\succeq \frac{1}{r \cdot e} \cdot \lambda$.
\end{proposition}

\begin{proof}
If $e = 1$, the sheaf $\sL |_{F'}$ is trivial, hence $f_* \omega^{\nu}_{X/Y}
= \lambda$ and \ref{bounds} obviously holds true. Hence we will
assume $e \geq 2$. For some $\mu \gg 0 $  there exists an
effective divisor $\Sigma_1$, disjoint from
$S$ with $\lambda^{\mu} = \sO_Y(\Sigma_1)$. By \ref{twist} and
by flat base change, we are free to replace $Y $ by any $Y'$,
finite over $Y$ and unramified over a neighborhood of $S$. Hence
we are allowed to assume that $\Sigma_1 = (\nu -1) \cdot e \cdot
\mu \cdot \Sigma$ or that
$$
\lambda = \sO_Y ((\nu - 1) \cdot e \cdot \Sigma).
$$
Consider the $r$-fold fibre product
$$
f^r : X^r = X'\times_Y X' \ldots \times_Y X' \>>> Y.
$$
$f^r$ is flat and Gorenstein and smooth over some open
subscheme. Let
$$\pi : X^{(r)} \>>> X^r$$
be a desingularization
such that the general fibre $F^{(r)}$ of $f^{(r)} = f^r \circ
\pi$ is isomorphic to $F \times \ldots \times F$. For
$$
\sM = \pi^* \bigotimes^{r}_{i=1} pr^{*}_{i} \sL
\subset \pi^*\omega^\nu_{X^r/Y},
$$
using flat base change, and the natural maps
$$
\sO_{X^r} \to \pi_*\sO_{X^{(r)}}\mbox{ \ \  and \ \ }
\pi_*\omega_{X^{(r)}} \to \omega_{X^{(r)}},
$$
one finds
\begin{gather}
\bigotimes^{r} f'_* \omega^{\nu}_{X'/Y} = \bigotimes^r \sE =
\bigotimes^r f_* \sL \>>> f^{(r)}_{*} \sM \label{nat-one} \ \
\mbox{ and }\\[.1cm]
f^{(r)}_{*} ((\pi^* \omega^{\nu-1}_{X^r /Y} )
\otimes \omega_{X^{(r)} /Y}) \>>> f^{r}_{*} \omega^{\nu}_{X^r/Y}
= \bigotimes f'_* \omega^{\nu}_{X'/Y}, \label{nat-two}
\end{gather}
and both are isomorphism over some open dense subset of $Y$.
(\ref{nat-one}) induces a surjection
$$
f^{(r)*} \bigotimes^r \sE = \pi^{*} \bigotimes^{r}_{i=1}
pr^{*}_{i} {f'}^* \sE \>>> \pi^* \bigotimes^{r}_{i=1}
pr^{*}_{i} \sL = \sM.
$$
In particular, since $\lambda \subset \otimes^r \sE$,
the sheaf ${f^{(r)}}^* \lambda$ is a subsheaf of $\sM$. Let
$\Gamma$ denote the divisor with $\sM (- \Gamma) =
f^{(r)*} \lambda$. For some divisor $C$, supported in fibres of
$f^{(r)}$ one has
$$\pi^* \omega_{X^r/Y} = \omega_{X^{(r)}/Y} (C).$$
Blowing up $X^{(r)}$ with centers in fibres of $f^{(r)}$ we find
a normal crossing divisor $D$ with
$$
D \geq \pi^* \sum^{r}_{i=1} pr^{*}_{i} B
$$
and such that $\sM (D) = \omega_{X^{(r)}/Y} (C)^{\nu}$.
For
$$\nabla = e \cdot (\nu -1) \cdot D + \nu \cdot \Gamma + e
\cdot \nu \cdot (\nu-1) \cdot f^{(r)*} (\Sigma)$$
one obtains
$$
\omega_{X^{(r)}/Y} (C)^{e \cdot \nu \cdot (\nu-1)} (- \nabla)
= \sM^{e \cdot (\nu -1)} (- \nu \cdot \Gamma) \otimes f^{(r)*}
\lambda^{-\nu} = \sM^{e \cdot(\nu -1) - \nu}.
$$
Since we assumed $e, \nu \geq 2$, the exponent of $\sM$ is
non-negative. By \ref{weak}, i), the sheaf
$\otimes^r\sE$ is nef and \ref{fujita} implies that
$$
\sF = f^{(r)}_{*} \left(\omega_{X^{(r)}/Y} (C)^{\nu-1} \otimes
\omega_{X^{(r)}/Y} \left\{ - \frac{\nabla}{e \cdot \nu} \right\}
\right)
$$
is nef. $\sF$ is contained in
$$
\sF' = f^{(r)}_{*} (\pi^* \omega^{\nu-1}_{X^r/Y} \otimes
\omega_{X^{(r)}/Y} ) \otimes \sO_Y (-(\nu-1) \cdot \Sigma),
$$
and using (\ref{nat-two}) ones finds
$$
\sF \subseteq (\bigotimes^r f'_* \omega^{\nu}_{X'/Y} ) \otimes \sO_Y (-(\nu-1)
\cdot \Sigma).
$$
On the other hand, $\sF$ contains
$$
\sF'' = f^{(r)}_{*} \left(\omega_{X^{(r)}/Y} (C)^{\nu -1} \otimes
\omega_{X^{(r)}/Y} \left\{ - \frac{\nabla - e \cdot D}{e \cdot
\nu} \right\} \right).
$$
Over some sufficiently small open dense subset $U \subset Y$
$$
\sF'' |_U = f^{(r)}_{*} \left(\omega^{\nu -1}_{X^{(r)}/Y} \otimes
\omega_{X^{(r)}/Y} \left.\left\{ - \frac{\nu \cdot \Gamma}{\nu \cdot
e} \right\} \otimes \sO_{X^{(r)}} (-D)\right)\right|_U.
$$
By definition $e = e (\sL |_F)$, and from \cite{E-V2} or \cite{Vie}, 5.21,
one has
$$
e = e (\sL |_F) = e (\bigotimes^{r}_{i=1} pr^{*}_{i} \sL |_F ) =
e (\sM |_{F^{(r)}} ).
$$
The semicontinuity of $e$ in \cite{Vie}, 5.14, implies that for $U$
small enough,
$$
\sF'' |_U = f^{(r)}_{*} (\omega^{\nu}_{X^{(r)}/Y} \otimes
\sO_{X^{(r)}} (-D)) = \sF' |_U .
$$
Hence $\sF, \sF'$ and $\sF''$ have the same rank and $\sF'$ is
nef. We obtain
$$
\bigotimes^r f'_* \omega^{\nu}_{X'/Y} \succeq \sO_Y ((\nu -1)
\cdot \Sigma ) = \frac{1}{e} \lambda
$$
or $f'_* \omega^{\nu}_{X'/Y} \succeq \frac{1}{r\cdot e}
\lambda$, as claimed.
\end{proof}

\section{Differential forms and cyclic coverings}

Let $Y$ be a non-singular projective curve of genus $g$, and let
$S \subset Y$ be a reduced divisor of degree $s$. We consider
again an $(n+1)$-dimensional manifold $X$ and a surjective
morphism $f: X \to Y$ with connected general fibre $F$. For
$\Delta = f^{-1} (S)$ we will assume that $f_0 = f |_{X -
\Delta}$ is smooth.

For some $\nu >1$, with $f_*\omega^\nu_{X/Y}\neq 0$,
we choose, as in \ref{constants}, a blowing
up $\tau : X' \to X$ and a normal crossing divisor $B$.
Hence writing $f' = f \circ \tau$ and $\sL = \omega^{\nu}_{X'/Y} (-
B)$ one has $f'_* \sL =f'_* \omega^{\nu}_{X'/Y}$ and $f{'}^*
f'_* \sL \twoheadrightarrow \sL$. Choose $S' \supset S$, such
that $f'$ is smooth outside of $\Delta' = f{'}^{-1} (S')$, and
such that $B - (B \cap \Delta')$ is a relative normal crossing divisor
over $Y - S'$. Blowing up with centers in $\Delta'$, we may also
assume $B + \Delta'$ to be a normal crossing divisor.

Let $\sA$ be an ample invertible sheaf with
\begin{equation} \label{ass-nef}
f'_* \omega^{\nu}_{X'/Y} \otimes \sA^{-\nu} = f_*
\omega^{\nu}_{X/Y} \otimes \sA^{-\nu} \mbox{ \ \ \ \ \ \ ample.}
\end{equation}
For $N = \nu \cdot \mu$ and $\mu$ sufficiently large, and for
$\sM = \omega_{X'/Y} \otimes f{'}^* \sA^{-1}$ the sheaf
$$\sM^N (- \mu \cdot B) = \sL^{\mu} \otimes f{'}^* \sA^{-N}$$
is generated by global sections. For a general section $\sigma$
with zero divisor $V(\sigma)$, the divisor $M = \mu B + V
(\sigma)$ as well as $M + \Delta'$ are normal crossing divisors.
Enlarging $S'$, if necessary, and replacing $X'$ by a blowing up
with centers in $\Delta'$, we may assume that all the fibres of
$$f\circ \tau: X'-\Delta'\>>> Y-S'
$$
intersect $M$ transversely.

As explained in \cite{E-V2}, \S \ 2, $\sigma$ defines a cyclic
covering $Z$ of $X'$. In explicit terms, writing $\sM^{(-i)} =
\sM^{-i} ([\frac{i \cdot M}{N}])$, the covering is given by
$$
\gamma : Z = {\rm \bf Spec} \left(\bigotimes^{N-1}_{i=0} \sM^{(-i)}\right)
\>>> X'.
$$
Let us write $g=f'\circ \gamma$ for the induced morphism from $Z$
to $Y$. The discriminant of $\gamma$ is $M' = (M - N \cdot [
\frac{M}{N}])_{{\rm red}}$, and defining
$$
\Omega^p_Z(\log \gamma^{-1}(\Delta'+M'))\mbox{ \ \  and \ \ }
\Omega^p_{Z/Y}(\log \gamma^{-1}(\Delta'+M'))
$$
to be the reflexive hulls of the
corresponding sheaves on the smooth locus of $Z$, both sheaves
are locally free and, as explained in \cite{E-V2}, \S 2, one has
\begin{equation}\begin{split}\label{d1}
\Omega^p_Z(\log \gamma^{-1}(\Delta'+M')) = {\gamma}^*\Omega^p_{X'}(\log
(\Delta'+M'))\\
\mbox{and \ \ \ }
\Omega^p_{Z/Y}(\log \gamma^{-1}(\Delta'+M')) = {\gamma}^*\Omega^p_{X'/Y}(\log
(\Delta'+M')).
\end{split}\end{equation}
Hence $\Omega^p_{X'}(\log (\Delta'+M'))\otimes \sM^{(-1)}$ is a direct
factor of $\gamma_*\Omega^p_Z(\log \gamma^{-1}(\Delta'+M'))$.
The image of the induced map
$$
\gamma^*\Omega^p_{X'}(\log \Delta'+M')\otimes \sM^{(-1)}\>>>
\Omega^p_Z(\log \gamma^{-1}(\Delta'+M'))
$$
lies in the subsheaf $\Omega^p_Z(\log \gamma^{-1}(\Delta'))\tilde{\ }$,
where $( \ \ )\tilde{\ }$ denotes the reflexive hull. Again, similar
inclusions hold true for the relative differential forms
over $Y$.

Kawamata's covering construction
(see \cite{Vie}, 2.2) allows to choose a finite covering $Z''\to Z$
with ${Z''}$ non singular. Blowing up centers in fibres over $Y$, we obtain
a non-singular variety $Z'$ and a generically finite map $\eta:
Z' \to Z$, such that all the fibres of
$$
g' = g\circ \eta : Z' \>>> Y
$$
are normal crossing divisors. By abuse of notations we
will add some points to $S'$, hence some fibres to $\Delta'$,
and assume that $g'$ is smooth outside of $\Pi' = g{'}^{-1}
(S')$ and that $\eta$ is finite outside of $\Pi'$. Let
$\gamma': Z' \to X'$ be the induced map. Since
$\gamma'$ factors through a non-singular variety, finite over
$X'$, the higher direct images $R^q\gamma'_*\sO_{Z'}$ are zero
for $q>0$.

Let ${M''}$ denote the proper transform of $M'$ in $Z'$.
Since $\gamma'$ is finite outside of $\Pi'$ one obtains
natural inclusions (see \cite{E-V2}, 3.20, for example)
\begin{equation}\begin{split}\label{d2}
\eta^* \Omega^{p}_{Z} (\log \gamma^{-1}(\Delta'+M')) \>>>
\Omega^{p}_{Z'} (\log (\Pi'+M''))
\mbox{ \ \ \ and}\\
\eta^* \Omega^{p}_{Z/Y} (\log \gamma^{-1}(\Delta'+M')) \>>>
\Omega^{p}_{Z'/Y} (\log (\Pi'+ M'')).
\end{split}\end{equation}
(\ref{d1}) and (\ref{d2}) together induce inclusions
\begin{equation}\begin{split}\label{d3}
{\gamma'}^*(\tau^*\Omega^p_{X/\bullet}(\log \Delta))\otimes\sM^{(-1)}\>>>
{\gamma'}^*\Omega^p_{X'/\bullet}(\log (\Delta'+M'))\otimes\sM^{(-1)}\\ \>>>
\eta^* \Omega^{p}_{Z/\bullet} (\log \gamma^{-1}(\Delta'+M')) \>>>
\Omega^{p}_{Z'/\bullet} (\log (\Pi'+{M''})),
\end{split}\end{equation}
where $\bullet$ stands for $Y$ or for ${\rm Spec}\C$, respectively.

Again, the image of composite of the injections in \ref{d3}
must have trivial residues along the components of ${M''}$,
hence one finds a natural map
$$
\iota_\bullet:{\gamma'}^*(\tau^*\Omega^p_{X/\bullet}(\log
\Delta))\otimes\sM^{(-1)}\>\subset >> \Omega^{p}_{Z'/\bullet} (\log \Pi').
$$
Consider the tautological sequence
\begin{equation}\begin{split}\label{b}
0 \>>> f^* \Omega^{1}_{Y} (\log S) \otimes \Omega^{p-1}_{X/Y}
(\log \Delta) \>>> \Omega^{p}_{X} (\log \Delta) \>>> \Omega^{p}_{X/Y} (\log
\Delta)\>>> 0.
\end{split}\end{equation}
Pulling it back to $X'$ and tensorizing with $\sM^{(-1)}$ one
obtains
\begin{equation}\begin{split}\label{c}
0 \>>> {f'}^* \Omega^{1}_{Y} (\log S) \otimes \tau^*\Omega^{p-1}_{X/Y}
(\log \Delta)\otimes \sM^{(-1)}\\
\>>> \tau^*\Omega^{p}_{X} (\log \Delta)\otimes \sM^{(-1)} \>>>
\tau^*\Omega^{p}_{X/Y} (\log \Delta)\otimes \sM^{(-1)}\>>> 0.
\end{split}\end{equation}
The inclusions $\iota_Y$, $\iota_{{\rm Spec}\C}$ and
$$
\iota: \Omega^{1}_{Y} (\log S) \hookrightarrow
\Omega^{1}_{Y}(\log S').
$$
induce a morphism from the pullback of this exact sequence to
\begin{equation}\label{a2}
0 \>>> g{'}^*
\Omega^{1}_{Y} (\log S') \otimes \Omega^{p-1}_{Z'/Y} (\log \Pi')
\>>> \Omega^{p}_{Z'} (\log \Pi') \>>> \Omega^{p}_{Z'/Y} (\log \Pi')
\>>> 0
\end{equation}
Let us define
\begin{gather*}
E^{p,q} = R^q g'_* \Omega^{p}_{Z'/Y} (\log \Pi') \\
\mbox{and \ \ \ \ \ }
F^{p,q} = R^q f'_* (( \tau^* \Omega^{p}_{X/Y} (\log \Delta))
\otimes \sM^{(-1)}).
\end{gather*}
The inclusion $\iota_Y$ gives a map
$$
R^qg'_*({\gamma'}^*(\tau^*\Omega^p_{X/Y}(\log \Delta)) \otimes
\sM^{(-1)}) \>>> R^qg'_*\Omega^{p}_{Z'/Y} (\log \Pi').
$$
Since the first sheaf is isomorphic to
$$
R^qf'_*(({\gamma'}_*{\gamma'}^*\sO_{Z'})\otimes
(\tau^*\Omega^p_{X/Y}(\log \Delta)) \otimes
\sM^{(-1)})
$$
we obtain thereby a morphism $\rho_{p,q}:F^{p,q} \to E^{p,q}$. Obviously
$$
\rho_{n,0}: f'_* (( \tau^* \Omega^{n}_{X/Y} (\log \Delta))
\otimes \sM^{(-1)}) \>>> g'_* \Omega^{n}_{Z'/Y} (\log \Pi')
$$
is injective, and
$$
\rho_{0,n}:R^n f'_* (\sM^{(-1)}) \>>> R^nf'_*(\gamma_*\sO_Z)
\>>> R^nf'_*(\gamma'_*\sO_{Z'})=R^n g'_* \sO_{Z'}
$$
gives $F^{0,n}$ as a direct factor of $E^{0,n}$.
The edge-morphism
$$
E^{p,q} \> \theta_{p,q} >> E^{p-1, q+1} \otimes \Omega^{1}_{Y}
(\log S')
$$
of the exact sequence (\ref{a2}) is the Kodaira Spencer map,
studied in \S 1. Since the pullback of (\ref{c}) to $Z'$ is a
subsequence of (\ref{a2}) $\theta_{p,q}$ commutes with the
edge-morphism
\begin{multline*}
R^qg'_*({\gamma'}^*(\tau^*\Omega^p_{X/Y}(\log \Delta)) \otimes
\sM^{(-1)})\\
\>>>
R^{q+1}g'_*({\gamma'}^*(\tau^*\Omega^{p-1}_{X/Y}(\log \Delta)) \otimes
\sM^{(-1)})\otimes \Omega^1_Y(\log S).
\end{multline*}
So the edge-morphism of the exact sequence (\ref{c}), denoted by
$$
F^{p,q} \> \tau_{p,q} >> F^{p-1,q+1} \otimes \Omega^{1}_{Y}
(\log S),
$$
is compatible with $\theta_{p,q}$.
\begin{lemma} \label{commutative} Assume for an ample
invertible sheaf $\sA$, and for $\nu>1$ (\ref{ass-nef}) holds true.
Then, using the notations introduced above,
\begin{enumerate}
\item[i)] the Kodaira Spencer map for $g' : Z' \to Y$ and the
edge-morphism of the exact sequence (\ref{c}) induce a
commutative diagram
$$
\begin{CD}
E^{p,q} \> \theta_{p,q} >> E^{p-1, q+1} \otimes \Omega^{1}_{Y}
(\log S') \\
\A \rho_{p,q} AA \A A \rho_{p-1,q+1} \otimes \iota A \\
F^{p,q} \> \tau_{p,q} >> F^{p-1,q+1} \otimes \Omega^{1}_{Y} (\log S).
\end{CD}
$$

\item[ii)] $\rho_{n,0}$ is injective.
\item[iii)] $F^{0,n}$ is a direct factor of $E^{0,n}$.
\item[iv)] the sheaf $(F^{0,n})^\vee$ is ample.
\item[v)] the sheaf $F^{n,0}$ is invertible of degree
$$\deg (F^{n,0}) \geq \deg (\sA) -\delta,$$
where $\delta$ denotes the number of non-semistable fibres.
\end{enumerate}
\end{lemma}

\begin{proof}
It remains to verify iv) and v). Comparing the first Chern
classes of the sheaves in (\ref{b}) one finds
$\Omega^{n}_{X/Y} (\log \Delta) = \omega_{X/Y} (\Delta_{{\rm
red}} - \Delta )$. Hence for some invertible sheaf $\sH$ of
degree $\delta$ one has an inclusion
$$\Omega^{n}_{X/Y} (\log \Delta) \supset \omega_{X/Y} \otimes f^*
\sH^{-1}.$$
Recall that $\sM = \omega_{X'/Y} \otimes f{'}^* \sA^{-1}$, and that
$$
\sM^{(-1)} = \sM^{-1} \left(\left[\frac{M}{N}\right]\right) =
\sM^{-1} \left( \left[\frac{B}{\nu} \right]\right),
$$
where $B$ is the relative fix locus of $\omega^{\nu}_{X'/Y}$.
In particular, if $E$ denotes the effective divisor with
$\omega_{X'/Y}(-E) = \tau^*\omega_{X/Y}$, the divisor $B$ is
larger than $\nu\cdot E$.

In order to prove that $\sA\otimes \sH^{-1}$
is a subsheaf of $F^{n,0}$, and that the latter is invertible,
we just have to show that
$$
f'_* (\tau^* \omega_{X/Y} \otimes \sM^{(-1)}) = f'_* \left((f{'}^{*}
\sA ) \otimes \sO_{X'} \left( \left[ \frac{B}{\nu} \right]
-E \right)\right)
$$
is isomorphic to $\sA$, or that $\sO_Y \simeq f'_* \sO_{X'}
\left( \left[ \frac{B}{\nu}\right] - E \right)$. This is
obvious since $0 \leq \left[ \frac{B}{\nu}\right] - E \leq B$
and since $B$ is the relative fix locus of an invertible sheaf.

\ref{kernel} says in particular, that $E^{0,n}={\rm Ker}(\tau_{0,n})$
has no invertible subsheaf of positive degree. Using iii) one
obtains the same for $F^{0,n} =R^n{f'}_*\sM^{(-1)}$
and $(F^{0,n})^\vee$
is nef. Serre duality and the projection formula imply
$$(F^{0,n})^\vee = \sA^{-1}\otimes {f'}_* (\omega_{X'/Y}^2\otimes
\sO_{X'}(-\big[\frac{B}{\nu}\big])).$$
To prove the ampleness, claimed in iv), choose some $\eta>0$
such that
$$
S^\eta(f_*\omega_{X/Y}^\nu\otimes \sA^{-\nu})\otimes \sA^{-1}
$$
is ample and consider a finite covering $\varphi:Y'\to Y$
of degree $\nu\cdot \eta$, \'etale over a neighborhood of the
discriminant divisor $S'$. Then $\varphi^*\sA =
{\sA'}^{\nu\cdot\eta}$, for some ample invertible sheaf $\sA'$
on $Y'$, and both, $f_*\omega^\nu_{X/Y}$ and $F^{0,n}$ are
compatible with pullbacks.

Replacing $Y$ by $Y'$, we may assume thereby that
$\sA = {\sA'}^{\nu\cdot\eta}$, for some $\sA'$, and that
$f_*\omega_{X/Y}^\nu\otimes \sA^{-\nu}\otimes {\sA'}^{-\nu}$
is ample. Repeating the argument used above, for
$\sA\otimes\sA'$ instead of $\sA$, one finds
$$\sA^{-1}\otimes{\sA'}^{-1}\otimes {f'}_* (\omega_{X'/Y}^2\otimes
\sO_{X'}(-\big[\frac{B}{\nu}\big]))$$
to be nef, hence $\sA^{-1}\otimes {f'}_* (\omega_{X'/Y}^2\otimes
\sO_{X'}(-\big[\frac{B}{\nu}\big]))$ to be ample.
\end{proof}

\section{The proof of \ref{main1}, \ref{main3} and \ref{main2}}

Recall that starting from $f: X \to Y$ we constructed a family
$g': Z' \to Y$, with discriminant locus $S' \supset S$. As in
\S 1, we consider the Hodge bundles
$$
E^{p,q} = R^q g'_* \Omega^{p}_{Z'/Y} (\log \Pi').
$$
In the last section we obtained a morphism
$$
\rho_{p,q}:F^{p,q} = R^q f'_* (\tau^*\Omega^{p}_{X/Y} (\log
\Delta) \otimes \sM^{(-1)}) \>>> E^{p,q},
$$
compatible with the Kodaira Spencer map
$$
\theta_{p,q}:E^{p,q} \to E^{p-1,q+1}\otimes \Omega^1_Y(\log
S').
$$
Moreover, by \ref{commutative} the image of $\rho_{p-1,q-1}\otimes \iota$
is contained in $E^{p-1,q+1}\otimes \Omega^1_Y(\log S)$.

\begin{proposition} \label{bounds-two}
Let $\delta$ denote the
number of those singular fibres of $f$ which are not reduced
normal crossing divisors and let $\nu > 1$ be an integer with
$f_* \omega^{\nu}_{X/Y} \neq 0$. If $\sA$ is an invertible
sheaf, with $\deg( \sA) > \delta$, and such that
$f_*\omega^{\nu}_{X/Y} \otimes \sA^{-\nu}$ is ample,
then
\begin{enumerate}
\item[a)] $(2g - 2 + s) \geq 0$ implies that $\deg (\sA) \leq n \cdot
(2 g - 2 + s ) +\delta$.
\item[b)] $(2g - 2 + s) > 0$ implies that $\deg (\sA) < n \cdot
(2 g - 2 + s ) +\delta$.
\end{enumerate}
\end{proposition}
\begin{proof} To handle both cases at once, define $\epsilon=1$,
if $2g-2+s=0$ and $\epsilon=0$, otherwise.
Assume that
$$\deg (\sA) \geq n \cdot (2 g - 2 + s ) + \epsilon
+\delta.
$$
Then by \ref{commutative}, v), $F^{n,0}$ is an
invertible subsheaf of $E^{n,0}$ of degree
$$
\deg (F^{n,0}) \geq n \cdot (2g - 2 + s) +\epsilon,
$$
in particular it is ample.
For $0 \leq i \leq n$, we will construct by induction
an invertible subsheaf $\tilde{F}^{n-i,i}$ of $\rho_{n-i,i}
(F^{n-i,i}) \subset E^{n-i,i}$ of degree
$$
\deg (\tilde{F}^{n-i,i}) \geq
(n-i) \cdot (2g - 2 + s) + \epsilon.
$$
If $i< n$, the sheaf $\tilde{F}^{n-i,i}$ is ample, and by \ref{kernel} it
can not lie in the kernel of $\theta_{n-i,i}$. On the
other hand, by \ref{commutative}, i),
$$
\theta_{n-i,i} (\tilde{F}^{n-i,i}) \subseteq \theta_{n-i,i}
(\rho_{n-i,i} F^{n-i,i} ) \subseteq \rho_{n-i-1,i+1} (F^{n-i-1,
i+1} )\otimes \Omega^{1}_{Y} (\log S).
$$
The invertible sheaf
$$
\tilde{F}^{n-i-1,i+1} = \theta_{n-i,i} (\tilde{F}^{n-i,i})
\otimes \Omega^{1}_{Y} (\log S)^{-1}
$$
thereby is a subsheaf of
$$
\rho_{n-i-1,i+1} (F^{n-i-1, i+1}) \subseteq E^{n-i-1,i+1}
$$
of degree $\deg (\tilde{F}^{n-i-1,i+1}) \geq (n - i -1) \cdot (2g - 2
+s)+\epsilon $.

For $i=n$ we obtain a subsheaf $\tilde{F}^{0,n}$ of degree
$\deg (\tilde{F}^{0,n}) \geq \epsilon \geq 0$,
contradicting \ref{commutative}, iv).
\end{proof}

Now everything is set to prove \ref{main1}, \ref{main2}, and
\ref{main3}. We will proceed in the following way: Adding one or
two points to $S$, if necessary, hence declaring some of the
smooth fibres to be ``singular'', we are allowed to assume
\begin{equation}
(2g - 2 + s) = \deg (\Omega^{1}_{Y} (\log S)) \geq 0.
\end{equation}
If $2g-2+s=0$ and if $f$ is semistable, one finds
the degree of $\lambda_{\nu}$ to be zero.
Then the assumptions a) and b) in \ref{main1} imply that
$f$ is isotrivial. For $Y=\P^1$, independently of the additional
assumptions on $F$, this will imply that $\kappa(X)=-\infty$.

For $2g-2+s>0$ the inequalities stated in
\ref{main2} i), ii), will hold true, independently
of the semi-ampleness of $\omega_F$, whenever $\deg(\lambda_\nu)
>0$.

\begin{proposition}\label{bounds-general1}
Let $Y_0$ be either an elliptic curve or $\C^*$, and let
$\tau:Y'\to Y$ be a finite morphism, \'etale over $Y_0$, such
that $X\times_YY' \to Y'$ has a semistable model $f':X'\to Y'$.
Then, for all $\nu >1$ with $H^0(F,\omega_F^\nu)\neq 0$, the degree
of $\det (f'_* \omega^{\nu}_{X'/Y'})$ is zero.
\end{proposition}

\begin{proof} The sheaf $\lambda':=
\det({f'}_*\omega_{X'/Y'}^\nu)$ is nef, hence if
\ref{bounds-general1} does not hold true,
it is ample. Since $f':X'\to Y'$ is semistable, $\lambda'$ is
compatible with further pullbacks. Replacing $Y'$ by a
covering, we may assume that
$\deg(\lambda') > \nu \cdot e \cdot r$. By Proposition \ref{bounds}
$$
f'_* \omega^{\nu}_{X'/Y'} \succeq \frac{1}{r \cdot e} \cdot \lambda',
$$
hence $f'_* \omega^{\nu}_{X'/Y'} \otimes \sO_{Y'}(-\nu\cdot y')$
is ample, for a point $y'\in Y'$.
\ref{bounds-two}, a) implies that
$\deg(\sO_{Y'}(y')) \leq 0$, obviously a contradiction.
\end{proof}

\begin{proposition} \label{bounds-general}
Assume that $f_*\omega^{\nu}_{X/Y} \neq 0$, and that
$\lambda = \det (f_* \omega^{\nu}_{X/Y} )$ is ample,
for some $\nu > 1$. Let $\delta$ denote the number of
non-semistable fibres, $r = \rank (f_* \omega^{\nu}_{X/Y})$,
and let $e$ be the constant introduced in \ref{constants}.
If $2g-2+s > 0$, then
$\deg (\lambda) \leq (n \cdot (2g - 2 +s) + \delta ) \cdot \nu
\cdot e \cdot r.$
\end{proposition}

\begin{proof}
Choose an invertible sheaf $\sA$ on $Y$ of
degree $n \cdot (2g - 2 +s ) + \delta$.
If
$$\deg (\lambda) > (n \cdot (2g - 2 + s) + \delta) \cdot \nu \cdot
e \cdot r,$$
one finds $\deg(\sA^{\nu \cdot e \cdot r}) < \deg(\lambda)$.
By Proposition \ref{bounds}
$$
f_* \omega^{\nu}_{X/Y} \succeq \frac{1}{r \cdot e} \cdot \lambda,
$$
hence $f_* \omega^{\nu}_{X/Y} \otimes \sA^{-\nu}$ is ample.
Proposition \ref{bounds-two}, b) implies that
$$
\deg (\sA) < n \cdot(2g - 2 + s) +\delta,
$$
contradicting the choice of $\sA$.
\end{proof}

\noindent  {\it Proof of} \ \ref{main1}.
Assume $f:X\to Y$ is a morphism, smooth over $Y_0$. If $Y_0$ is
an elliptic curve, $f$ is smooth. For $Y_0=\C^*$, there exists
a finite covering $\tau:Y'\to Y$, \'etale over $\C^*$,
and a semistable family $f':X'\to Y'$, birational to the
pullback of $f$. Hence to show that $f$ is isotrivial,
\ref{bounds-general1} allows to assume that
$\deg(\det (f_*\omega^{\nu}_{X/Y}))=0$, for all $\nu>1$.

The experts will have noticed that the assumption a) or b) in
\ref{main1} are exactly those needed by Kawamata, Koll\'ar and
the first named author to prove the additivity of the Kodaira
dimension, and even the stronger statement $Q_{n,1}$, saying that
for non-isotrivial morphisms $f$
$$
\kappa (\omega_{X/Y}) = \kappa( F) + 1.
$$
Using \ref{weak}, this is equivalent to the ampleness of
$\det (f_*\omega^{\nu}_{X/Y} )$, for all positive multiples
$\nu$ of some $\nu_0 \gg 1$ (see \cite{Mor}, 7.2 and 7.6, and the
references given there).
In particular, the morphism $f:X\to Y$, considered above, is
isotrivial. \qed \\

\noindent {\it Proof of} \ \ref{main3}.
Assume there exist a morphism $f: X \to \P^1$, smooth over $\C^*$,
and with $\kappa(X) \geq 0$.
Let
$$
X\> f' >> Y' \> \tau >> \P^1
$$
be the Stein factorization. $\tau$ must be smooth over
$\C^*$, hence $Y'=\P^1$ and $\tau^{-1}(S)=\{0,\infty\}$.

Altogether we find a morphism, denoted again by $f:X\to \P^1$,
which is smooth over $\C^*$ and whose general fibre $F$ is connected.

For a finite morphism $\P^1 \to \P^1$ of degree $d$, \'etale over
$\C^*$, let
$$\varphi:Z\to X\times_{\P^1}\P^1$$
be a desingularization, and $g=pr_2\circ\varphi:Z \to \P^1$. Then
$\varphi^*pr_1^*\omega^\nu_X$ is a subsheaf of $\omega_Z^\nu$,
hence $\kappa(Z) \geq \kappa(X) \geq 0$. In particular, for some
$\nu_0$ and for all positive multiples $\nu$ of $\nu_0$ the sheaf
$$ g_*\omega_Z^\nu= \sO_{\P^1}(-2\cdot\nu)\otimes
g_*\omega^\nu_{Z/\P^1} $$ has a non trivial section. Therefore
$g_*\omega^\nu_{Z/\P^1}$ contains a non-trivial ample subsheaf. On
the other hand, $g_*\omega^\nu_{Z/\P^1}$ is nef, hence
$\det(g_*\omega^\nu_{Z/\P^1})$ must be ample. Since $g$ is smooth
over $\P^1-\{0,\infty\}$, and semistable for $d$ sufficiently
large and divisible, this contradicts \ref{bounds-general1}.
\qed\\

\noindent {\it Proof of} \ \ref{main2}. In fact, by
\ref{bounds-general} the inequalities
i) and ii) hold true under the assumptions a) and b) stated in
\ref{main1}. However, the constant $e$, defined in
\ref{constants} might depend on the general fibre of $f$ and not just
on the Hilbert polynomial $h(t)$. For this reason we
have to require in \ref{main2} the sheaf $\omega_F$ to be
semi-ample.

Under this assumption, there exists a quasi-projective moduli
scheme $M_h$, parameterizing pairs $(F, \sH)$ where $F$ is a
manifold with $\omega_F$ semi-ample and with $\sH$ a
polarization with Hilbert polynomial $h$ (see \cite{Vie}).
Seshadri and Koll\'ar constructed a finite covering $Z$ of
$M_h$, which carries a universal family (see \cite{Vie}, 9.25),
containing all $(F, \sH)$ with Hilbert polynomial $h(t)$ as
fibres. So we find some $\nu > 0$, such that for all such $(F,
\sH)$, $\omega^{\nu}_{F}$ is generated by global sections, hence
we may choose the sheaf $\sL$ in \ref{constants} in such a way
that $\sL |_F = \omega^{\nu}_{F})$. Finally, by
\cite{Vie}, 5.17, $e(\sL|_F)= e(\omega_F^\nu)$ is upper
semicontinous for the Zariski topology. Hence there exists some
$e$, with $e \geq e(\omega_F^\nu)$ for all
$F$. Altogether, we can choose $\nu$ and $e$ in \ref{main2}, to
depend just on $M_{h}$, hence just on $h (t)$.\qed \\

\bibliographystyle{plain}

\end{document}